\documentclass{amsart}
\usepackage{graphicx}

\newtheorem{thm}{Theorem}[section]
\newtheorem{prp}{Proposition}[section]

\newtheorem{cor}{Corollary}[section]

\newtheorem{defn}{Definition}[section]
\newtheorem{obs}{Observation}[section]
\DeclareMathOperator{\Des}{Des}
\DeclareMathOperator{\des}{des}
\DeclareMathOperator{\cDes}{cDes}
\DeclareMathOperator{\cdes}{cdes}
\DeclareMathOperator{\aDes}{aDes}
\DeclareMathOperator{\ades}{ades}
\DeclareMathOperator{\sgn}{sgn}
\DeclareMathOperator{\maj}{maj}
\DeclareMathOperator{\comaj}{comaj}
\DeclareMathOperator{\acomaj}{acomaj}

\title{Cyclic Descents and $P$-partitions}

\author[T. K. Petersen]{T. Kyle Petersen}
\address{Department of Mathematics, Brandeis University, Waltham, MA, USA, 02454}
\email{tkpeters@brandeis.edu}
\urladdr{http://people.brandeis.edu/\~{}tkpeters}

\begin{document}
\begin{abstract}
Louis Solomon showed that the group algebra of the symmetric group
$\mathfrak{S}_{n}$ has a subalgebra called the descent algebra,
generated by sums of permutations with a given descent set. In
fact, he showed that every Coxeter group has something that can
be called a descent algebra. There is also a commutative, semisimple subalgebra of Solomon's descent algebra generated by sums of permutations with the same number of descents: an ``Eulerian" descent algebra. For any Coxeter group that is also a Weyl group, Paola Cellini proved the existence of a different Eulerian subalgebra based on a modified definition of descent. We derive the existence of Cellini's subalgebra for the case of the
symmetric group and of the hyperoctahedral group using a variation
on Richard Stanley's theory of $P$-partitions.
\end{abstract}

\maketitle

\section{Introduction}
\subsection{The symmetric group} \label{sec:sym}
The study of permutations by descent sets is a natural
generalization of the study of permutations by number of descents:
the study of Eulerian numbers. Let $\mathfrak{S}_{n}$ be the
symmetric group on $n$ elements. We think of permutations in $\mathfrak{S}_n$ as bijections \[\pi: [n] \to [n],\] where $[n]$ denotes the set $\{1,2,\ldots,n\}$. For any permutation $\pi \in
\mathfrak{S}_{n}$, we say $\pi$ has a \emph{descent} in position
$i$ if $\pi(i) > \pi(i+1)$. Define the set $\Des(\pi) = \{\, i
\mid 1\leq i\leq n-1, \pi(i) > \pi(i+1)\,\}$ and let $\des(\pi)$
denote the number of elements in $\Des(\pi)$. We call $\Des(\pi)$
the \emph{descent set} of $\pi$, and $\des(\pi)$ the \emph{descent
number} of $\pi$. For example, the permutation $\pi =
(\pi(1),\pi(2),\pi(3),\pi(4)) = (1,4,3,2)$ has descent set
$\{2,3\}$ and descent number 2. The number of permutations of $[n]$
with descent number $k$ is denoted by the Eulerian number
$A_{n,k+1}$, and we recall that the Eulerian polynomial is defined
as \[A_{n}(t) = \sum_{\pi\in\mathfrak{S}_{n}} t^{\des(\pi)+1} =
\sum_{i=1}^{n} A_{n,i}t^{i}.\]

For each subset $I$ of $[n-1]$, let $u_{I}$ denote
the sum, in the group algebra $\mathbb{Q}[\mathfrak{S}_{n}]$, of all
permutations with descent set $I$. Louis Solomon \cite{Sol} showed
that the linear span of the $u_{I}$ forms a subalgebra of the
group algebra, called the \emph{descent algebra}. More generally,
he showed that one can define a descent algebra for any Coxeter
group.

For now consider the descent algebra of the symmetric group. This
descent algebra was studied in great detail by Adriano Garsia and Christophe
Reutenauer \cite{GR}. Jean-Louis Loday \cite{JLL} proved the existence of a commutative, semisimple subalgebra of Solomon's descent algebra. Sometimes called the ``Eulerian subalgebra,"
it is defined as follows. For $1\leq i\leq n$, let $E_{i}$ be the
sum of all permutations in $\mathfrak{S}_{n}$ with descent number
$i-1$. Then the Eulerian descent algebra is the linear span of the $E_i$. Define
\[ \phi(x) = \sum_{\pi\in\mathfrak{S}_{n}}
\binom{x+n-1-\des(\pi)}{n}\pi =
\sum_{i=1}^{n}\binom{x+n-i}{n}E_{i},\] which we refer to as the ``structure polynomial." Notice that the structure polynomial is a polynomial in $x$, with coefficients in $\mathbb{Q}[\mathfrak{S}_n]$, of degree $n$ and with no constant term. (In other words, it has exactly as many nonzero terms as there are possible descent numbers.) The Eulerian subalgebra is described by the following multiplication of structure polynomials:
\begin{thm}[Gessel]\label{thm:ges}
As polynomials in $x$ and $y$ with coefficients in the group algebra, we have
\begin{equation}\label{eq1}
\phi(x)\phi(y) =\phi(xy).
\end{equation}
\end{thm}
Define elements $e_{i}$ in the
group algebra by $\displaystyle \phi(x) = \sum_{i=1}^{n}
e_{i}x^{i}$. By examining the coefficients of $x^i y^j$ in \eqref{eq1}, it is clear that the $e_{i}$ are
orthogonal idempotents: $e_{i}^{2} = e_{i}$ and $e_{i}e_{j} = 0$
if $i\neq j$. This formulation is essentially an unsigned version of Loday's Th\'{e}or\`{e}me 1.7 (see \cite{JLL}). In particular, for $1 \leq k \leq n$, we have $\phi(k) = | \lambda^k |$, where Loday defines
\begin{eqnarray*}
  \lambda^k  & = & \sum_{ i = 0 }^{k-1} (-1)^i \binom{n+i}{i} l^{k-i}, \\
 l^{k} & = & (-1)^{k-1} \sum_{ \des{\pi} = k-1 } \sgn(\pi) \pi,
\end{eqnarray*}
and $| \lambda^k |$ means that we forget the signs on all the summands.

As an interesting aside, we remark that Theorem \ref{thm:ges} has connections to the well-known card shuffling analysis of Dave Bayer and Persi Diaconis \cite{BayDia}. If we define \[\overline{\phi}(x) = \sum_{\pi \in \mathfrak{S}_n} \binom{x+n-1-\des(\pi)}{n}\pi^{-1},\] then $\overline{\phi}(2)$ is the generating function for the probability distribution for one shuffle of a deck of $n$ cards (as defined in \cite{BayDia}).  We use the multiplication rule given by Theorem \ref{thm:ges} to compute the probability distribution after $m$ (independent) shuffles of the deck (it is not hard to show that $\overline{\phi}(x)$ obeys the same multiplication rule as $\phi(x)$): \[ \underbrace{ \overline{\phi}(2)\overline{\phi}(2)\cdots \overline{\phi}(2) }_{m \mbox{ \tiny   times}} = \overline{\phi}(2^m), \] which gives Theorem 1 from \cite{BayDia}.

Theorem \ref{thm:ges} can be proved in several ways, but we will focus on one that employs Richard Stanley's theory of $P$-partitions. More specifically, the approach taken in this
paper follows from work of Ira Gessel---Theorem \ref{thm:ges} is in
fact an easy corollary of Theorem 11 from \cite{G2}. Later we will
give a proof of Theorem \ref{thm:ges} that derives from Gessel's
work and then extend this method to prove the main results of this
paper. While equivalent results can be found elsewhere in the literature, the $P$-partition approach is self-contained and rather elementary. Further, it often yields $q$-analogs for formulas like \eqref{eq1} quite naturally.

Now we introduce another descent algebra based on a modified definition of descent. For a permutation $\pi \in \mathfrak{S}_{n}$ we define a \emph{cyclic descent} at position $i$ if $\pi(i) > \pi(i+1)$, or
if $i = n$ and $\pi(n) > \pi(1)$. Define $\cDes(\pi)$ to be the
set of cyclic descent positions of $\pi$, called the \emph{cyclic
descent set}. Let the \emph{cyclic descent number}, $\cdes(\pi)$,
be the number of cyclic descents. The number of cyclic descents is between 1 and $n-1$. One can quickly observe that a permutation $\pi$ has the same number of cyclic descents as both $\pi\omega^{i}$ and $\omega^i \pi$ for $i=0,1,\ldots,n-1$, where $\omega$ is the $n$-cycle $(1\,\,2\,\,\cdots\,n)$. Define the \emph{cyclic Eulerian polynomial} to be \[ A_{n}^{(c)}(t) = \sum_{\pi\in\mathfrak{S}_{n}}t^{\cdes(\pi)} = \sum_{i =1}^{n-1} A^{(c)}_{n,i}t^{i},\] where $A^{(c)}_{n,k}$ is the number of permutations with cyclic descent number $k$. We can make the following proposition (observed by Jason Fulman in \cite{F2}, Corollary 1).

\begin{prp}\label{prp:eul}
The cyclic Eulerian polynomial is expressible in terms of the ordinary Eulerian polynomial:
\[A_{n}^{(c)}(t) = nA_{n-1}(t).\]
\end{prp}

\begin{proof}
We will compare the coefficient of $t^d$ on each side of the equation to show $A^{(c)}_{n,d} = nA_{n-1,d}$. Let $\pi \in \mathfrak{S}_{n-1}$ be any permutation of $[n-1]$ such that $\des(\pi)+1 = d$. Let $\widetilde{\pi} \in \mathfrak{S}_{n}$ be the permutation defined by $\widetilde{\pi}(i) = \pi(i)$ for $i = 1,2,\ldots,n-1$ and $\widetilde{\pi}(n) = n$. Then we have $\des(\widetilde{\pi})= \des(\pi)$ and $\cdes(\widetilde{\pi}) = d$. Let $\langle\widetilde{\pi}\rangle = \{\, \widetilde{\pi}\omega^{i} \mid  i = 0,1,\ldots,n-1\,\}$, the set consisting of all $n$ cyclic permutations of $\widetilde{\pi}$. Every permutation in the set has exactly $d$ cyclic descents. There is a bijection between permutations of $\mathfrak{S}_{n-1}$ and such subsets of $\mathfrak{S}_n$ given by the map \[\pi \mapsto \langle\widetilde{\pi}\rangle, \] and so the proposition follows.
\end{proof}

Let $E_{i}^{(c)}$ be the sum in the group algebra
of all those permutations with $i$ cyclic descents. Then we define the cyclic structure polynomial
\[\varphi(x) =
\frac{1}{n}\sum_{\pi\in\mathfrak{S}_{n}}\binom{x+n-1-\cdes(\pi)}{n-1}\pi =
\frac{1}{n}\sum_{i=1}^{n-1}\binom{x+n-1-i}{n-1}E_{i}^{(c)}.\] Notice that $\varphi(x)$ is a polynomial of degree $n-1$ with no constant term, giving as many nonzero terms as possible cyclic descent numbers. Paola Cellini studied cyclic descents more generally in the papers \cite{Ce}, \cite{Ceii}, and \cite{Ce2}. Though her approach
is quite different from the one taken in this paper, from her work and Loday's we can derive the following theorem.
\begin{thm}\label{thmcyc}
As polynomials in $x$ and $y$ with coefficients in the group
algebra of the symmetric group, we have
\[\varphi(x)\varphi(y) = \varphi(xy).\]
\end{thm}
Now if we define elements $e_{i}^{(c)}$ by
$\displaystyle\varphi(x) = \sum_{i=1}^{n-1}e_{i}^{(c)}x^{i}$, we
see that $\left(e_{i}^{(c)}\right)^2 = e_{i}^{(c)}$ and
$e_{i}^{(c)} e_{j}^{(c)} = 0$ if $i\neq j$. Therefore the elements
$e_{i}^{(c)}$ are orthogonal idempotents. Fulman relates cyclic descents to card shuffling (now we consider cuts as well as shuffles of the deck); Theorem 2 of \cite{F2} is closely related to Theorem \ref{thmcyc}. We will prove
Theorem \ref{thmcyc} in section \ref{sec:cycprf} using formula \eqref{eq1}.

That the multiplication for the cyclic structure polynomials is so similar to that of the ordinary structure polynomials is no accident. Indeed, it will be clear from the proof of Theorem \ref{thmcyc} that the map \[ \pi \mapsto \sum_{ \sigma \in \langle\widetilde{\pi}\rangle } \sigma, \] where $\langle\widetilde{\pi}\rangle$ is as in the proof of Proposition \ref{prp:eul}, gives an isomorphism of algebras between the ordinary Eulerian descent algebra of $\mathbb{Q}[ \mathfrak{S}_{n-1} ]$ and the cyclic Eulerian descent algebra of $\mathbb{Q}[\mathfrak{S}_n ]$.

\subsection{The hyperoctahedral group} Let $\pm[n]$ denote the set
$\{-n,-n+1,\ldots,-1,0,1,\ldots, n-1,n\}$. Let $\mathfrak{B}_{n}$
denote the hyperoctahedral group, the group of all bijections
$\pi: \pm[n] \to \pm[n]$ with the property that $\pi(-s) =
-\pi(s)$, for $s = 0,1,\ldots,n$. Since the elements of the
hyperoctahedral group are uniquely determined by where they map
$1,2,\ldots,n$, we can think of them as signed permutations. For a
signed permutation $\pi \in \mathfrak{B}_{n}$ we will write $\pi =
(\pi(1),\pi(2),\ldots,\pi(n))$.

In moving from the symmetric group to the hyperoctahedral group, we define the \emph{descent set}
$\Des(\pi)$ of a signed permutation $\pi\in \mathfrak{B}_{n}$ to
be the set of all $i \in \{0,1,2,\ldots,n-1\}$ such that $\pi(i)
> \pi(i+1)$, where we always take $\pi(0) = 0$. The \emph{descent number} of
$\pi$ is again denoted $\des(\pi)$ and is equal to the cardinality of
$\Des(\pi)$.\footnote{It will be clear from the context whether we
are referring to the descent set of an ordinary permutation or
that of a signed permutation.} As a simple example, the signed
permutation $(-2,1)$ has descent set $\{0\}$ and descent number
$1$.

There is an Eulerian subalgebra for the hyperoctahedral group, previously studied by Cellini \cite{Ceii}, \cite{Ce2}, Francois and Nantel Bergeron \cite{BB1}, \cite{BB}, N. Bergeron \cite{B}, and Fulman \cite{F}. For $1\leq i \leq n+1$ let $E_{i}$ be the sum of all permutations in $\mathfrak{B}_{n}$ with $i-1$ descents. Define the type B structure polynomial
\[\phi(x) =
\sum_{\pi\in\mathfrak{B}_{n}}\binom{x+n-\des(\pi)}{n}\pi =
\sum_{i=1}^{n+1} \binom{x+n+1-i}{n} E_{i}.\] Chak-On Chow
\cite{Ch} was able to use the theory of $P$-partitions to prove the following.
\begin{thm}[Chow]\label{thm:cho}
As polynomials in $x$ and $y$ with coefficients in the group algebra of $\mathfrak{B}_n$, we have
\[\phi(x)\phi(y) = \phi(2xy+x+y),\]
or upon substituting \[x \leftarrow (x-1)/2,\] \[y \leftarrow (y-1)/2,\] then
\[
\phi((x-1)/2)\phi((y-1)/2) = \phi((xy-1)/2).
\]
\end{thm}
We therefore have orthogonal idempotents $e_{i}$ defined by $\displaystyle
\phi((x-1)/2) = \sum_{i=0}^{n} e_{i}x^{i}$ (notice that shifting the polynomial by 1/2 gives a nontrivial constant term). Just as Theorem \ref{thm:ges} was equivalent to Th\'{e}or\`{e}me 1.7 of \cite{JLL}, Theorem \ref{thm:cho} is equivalent to Theorem 4.1 from \cite{BB}. Further, both Fulman \cite{F} and Bergeron and Bergeron \cite{BB} made connections to ``signed" card shuffling (flip a card face up or down to change sign).

For a permutation $\pi \in \mathfrak{B}_{n}$, we define a
\emph{type B cyclic descent}, or \emph{augmented descent} at position $i$ if $\pi(i) > \pi(i+1)$ or
if $i=n$ and $\pi(n) > 0 = \pi(0)$. If we consider that signed
permutations always begin with 0, then augmented descents are the
natural generalization of type A cyclic descents.\footnote{Most
generally, Cellini \cite{Ce} uses the term ``descent in zero" to
represent this concept for any Weyl group. The term ``augmented descent" was coined by Gessel, \cite{G1}. Conceptually, the name ``type B cyclic descent" makes more sense, but in this paper we use the ``augmented" label, primarily because it requires fewer characters to typeset.} The set of all
augmented descent positions is denoted $\aDes(\pi)$, the
\emph{augmented descent set}. It is the ordinary descent set of
$\pi$ along with $n$ if $\pi(n) > 0$. The \emph{augmented descent
number}, $\ades(\pi)$, is the number of augmented descents. With
these definitions, $(-2,1)$ has augmented descent set $\{0,2\}$
and augmented descent number $2$. Note that while $\aDes(\pi)
\subset \{0,1,\ldots,n\}$, $\aDes(\pi) \neq \emptyset$, and
$\aDes(\pi) \neq \{0,1,\ldots,n\}$. Denote the number of signed
permutations in $\mathfrak{B}_n$ with $k$ augmented descents by $A_{n,k}^{(a)}$ and
define the \emph{augmented Eulerian polynomial} as
\[A_{n}^{(a)}(t) = \sum_{\pi\in\mathfrak{B}_{n}} t^{\ades(\pi)}
= \sum_{i=1}^{n}A_{n,i}^{(a)}t^{i}.\] In section \ref{sec:def} we
prove the following relation using the theory of $P$-partitions.

\begin{prp}\label{prp:augeul}
The number of signed permutations with $i+1$ augmented descents is
$2^{n}$ times the number of unsigned permutations with $i$
descents, $0\leq i \leq n-1$. In other words, $A_{n}^{(a)}(t) = 2^{n}A_{n}(t)$.
\end{prp}

Define the type B cyclic structure polynomial as
\[\psi(x) =
\sum_{\pi\in\mathfrak{B}_{n}}\binom{x+n-\ades(\pi)}{n}\pi = \sum_{i=1}^{n} \binom{x+n-i}{n} E^{(a)}_i,\] where $E^{(a)}_i$ is the sum of all signed permutations with $i$ augmented descents.

\begin{thm}\label{conj1}
As polynomials in $x$ and $y$ with coefficients in the group
algebra of the hyperoctahedral group we have
\[\psi(x)\psi(y) = \psi(2xy),\] or upon substituting \[x \leftarrow x/2,\] \[y \leftarrow y/2,\]
then \[\psi(x/2)\psi(y/2) = \psi(xy/2).\]
\end{thm}
We get orthogonal idempotents $e_{i}^{(a)}$ defined by
$\psi(x/2) = \displaystyle\sum_{i=1}^{n}e_{i}^{(a)}x^{i}$. We
will prove Theorem \ref{conj1} in section \ref{sec:prf} using
modified types of $P$-partitions called \emph{augmented}
$P$-partitions. Also in section \ref{sec:prf} we prove the
following theorem for which we have no type A equivalent.
\begin{thm}\label{thm:hypermod}
As polynomials in $x$ and $y$ with coefficients in the group algebra of the hyperoctahedral group we have
\[
\psi(x)\phi(y) = \psi(2xy+x),
\] or upon substituting \[x \leftarrow x/2,\] \[y \leftarrow (y-1)/2,\]
then \[\psi(x/2)\phi((y-1)/2) = \psi(xy/2).\]
\end{thm}
This formula tells us that $e_{i}^{(a)} e_i = e_{i}^{(a)}$ and that $e_{i}^{(a)} e_j = 0$ if $i\neq j$. In other words, the augmented Eulerian descent algebra is an ideal in the subalgebra of the group algebra formed by the span of the $e_i$ and the $e_{i}^{(a)}$. This relationship shows up again in the case of peak algebras of type A, and forms the basis of future work as discussed in section \ref{sec:fut}. See the paper of Marcelo Aguiar, N. Bergeron, and Kathryn Nyman \cite{AguBerNym} for more. We also point out here that when taken together, Theorems \ref{thm:cho}, \ref{conj1}, and \ref{thm:hypermod}, imply Cellini's Theorem A from \cite{Ceii}. Indeed, for positive integers $k$, Cellini's elements $x_{2k+1}$ are equivalent to $\overline{\phi}(k)$, and the $x_{2k+2}$ are equivalent to $\overline{\psi}(k)$. Here, as in section \ref{sec:sym}, $\overline{\phi}$ and $\overline{\psi}$ are the same as $\phi$ and $\psi$ except that we replace the coefficient of $\pi$ with the coefficient of $\pi^{-1}$. See also Lemma 4 of \cite{F}.

\section{The cyclic descent algebra}\label{sec:cyc}
\subsection{Definitions}\label{sec:def}
Throughout this paper we will use $P$ to denote a partially ordered set,
or poset. The partial order on $P$ is denoted $<_{P}$, or simply
$<$ when the meaning is clear. Our posets will always be finite
and for a poset of $n$ elements, the elements will be
distinctly labeled by the numbers $1,2,\ldots,n$.

\begin{defn}\label{def:p}
Let $X = \{x_{1},x_{2},\ldots \}$ be a countable, totally ordered
set. For a given poset $P$, a
\emph{$P$-partition}\footnote{We note that this definition varies from Richard Stanley's \cite{St} in that our maps are order-preserving, whereas his are order reversing, i.e., $f(i) \geq f(j)$ if $i <_P j$. However, in \cite{St} the maps of Definition \ref{def:p} are called (perhaps misleadingly) \emph{reverse} $P$-partitions.} is a function $f: [n] \to
X$ such that:
\begin{itemize}
\item $f(i) \leq f(j)$ if $i <_{P} j$

\item $f(i) < f(j)$ if $i <_{P} j$ and $i > j$ in $\mathbb{Z}$
\end{itemize}
\end{defn}
For our purposes we usually think of $X$ as a subset of the
positive integers. Let $\mathcal{A}(P)$ denote the set of all
$P$-partitions. When $X$ has finite cardinality $k$, then the number of
$P$-partitions is finite. In this case, define the \emph{order
polynomial}, denoted $\Omega_{P}(k)$, to be the number of
$P$-partitions $f:[n] \to X$.

We will consider any permutation $\pi \in \mathfrak{S}_{n}$ to be
a poset with the total order $\pi(s) <_{\pi} \pi(s+1)$. For
example, the permutation $\pi = (3,2,1,4)$ has $3 <_{\pi} 2
<_{\pi} 1 <_{\pi} 4$ as a poset. With this convention, the set
of all $\pi$-partitions is easily characterized. The set
$\mathcal{A}(\pi)$ is the set of all functions $f: [n] \to X$ such
that (if we take $X$ to be the positive integers)
\[1 \leq f(\pi_{1}) \leq f(\pi_{2})\leq \cdots \leq f(\pi_{n}),\]
and whenever $\pi(s) > \pi(s+1)$, then $f(\pi(s))< f(\pi(s+1))$,
$s=1,2,\ldots,n-1$. The set of all $\pi$-partitions where $\pi =
(3,2,1,4)$ is all maps $f$ such that $1 \leq f(3) < f(2) < f(1)
\leq f(4)$.

\begin{figure} [h]
\centering
\includegraphics{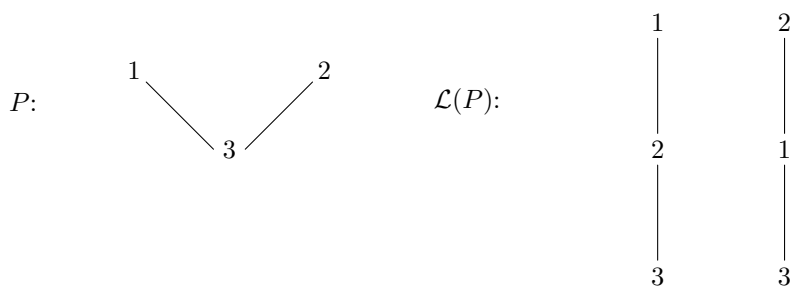}
\caption{Linear extensions of a poset $P$.}
\end{figure}

For a poset $P$ of order $n$, let $\mathcal{L}(P)$ denote the set of all
permutations of $[n]$ which extend $P$ to a total order. This set is sometimes called the set
of ``linear extensions" of $P$. For example let $P$ be the poset
defined by $1 >_{P} 3 <_{P} 2$. In ``linearizing" $P$ we form a
total order by retaining all the relations of $P$ but introducing
new relations so that any element is comparable to any other. In
this case, 1 and 2 are not comparable, so we have exactly two ways
of linearizing $P$: $3 < 2 < 1$ and $3 < 1 < 2$. These correspond
to the permutations $(3,2,1)$ and $(3,1,2)$. Let us make the
following observation.
\begin{obs}\label{ob1} A permutation $\pi$ is in $\mathcal{L}(P)$ if and only if $i<_{P}
j$ implies $\pi^{-1}(i) < \pi^{-1}(j)$.
\end{obs}
In other words, if $i$ is ``below" $j$ in the Hasse diagram of the poset $P$, it had better be below $j$ in any linear extension of the poset.
We also now prove what is sometimes called the fundamental theorem
(or lemma) of $P$-partitions.
\begin{thm}[FTPP]
The set of all $P$-partitions of a poset $P$ is the disjoint union
of the set of $\pi$-partitions of all linear extensions $\pi$ of
$P$: \[\mathcal{A}(P) = \coprod_{\pi \in \mathcal{L}(P)}
\mathcal{A}(\pi).\]
\end{thm}
\begin{proof}
The proof follows from induction on the number of incomparable
pairs of elements of $P$. If there are no incomparable pairs, then
$P$ has a total order and already represents a permutation.
Suppose $i$ and $j$ are incomparable in $P$. Let $P_{ij}$ be the
poset formed from $P$ by introducing the relation $i < j$. Then it
is clear that $\mathcal{A}(P) = \mathcal{A}(P_{ij}) \coprod
\mathcal{A}(P_{ji})$. We continue to split these posets (each with strictly fewer incomparable pairs)
until we have a collection of totally ordered chains corresponding to distinct linear
extensions of $P$.
\end{proof}
\begin{cor}
\[\Omega_{P}(k) = \sum_{\pi \in \mathcal{L}(P)} \Omega_{\pi}(k).\]
\end{cor}

The fundamental theorem tells us that in order to study $P$-partitions for a given poset, we can focus on the $\pi$-partitions for its linear extensions---a more straightforward task. In particular, counting $\pi$-partitions is not too difficult. Notice that for any permutation $\pi$ and any positive integer $k$,
\[\binom{k+n-1-\des(\pi)}{n} = \bigg(\!\!\binom{k-\des(\pi)}{n}\!\!\bigg),\] where
$\big(\!\binom{a}{b}\!\big)$ denotes the ``multi-choose" function---the number of ways to choose $b$ objects from a set of $a$ objects with repetitions. Another interpretation of $\big(\!\binom{a}{b}\!\big)$
is the number of integer solutions to the set of inequalities \[ 1\leq i_{1} \leq i_2 \leq \cdots \leq i_b \leq a. \] With this in mind, $\binom{k+n-1-\des(\pi)}{n}$ is the same as
the number of solutions to \[ 1\leq i_1 \leq i_2 \leq \cdots \leq
i_n \leq k-\des(\pi).\] Better still, we can say it is the number
of solutions (though not in general the same \emph{set} of solutions) to
\begin{equation}\label{eq:is}
1\leq i_1 \leq i_2 \leq \cdots \leq i_n \leq k \mbox{ \quad and  } i_s < i_{s+1} \mbox{ if } s \in \Des(\pi).
\end{equation}
For example, the number of solutions to $1 \leq i < j \leq 4$ is the same as the number of solutions to $1 \leq i \leq j-1 \leq 4-1$ or the solutions to $1 \leq i \leq j' \leq 3$. In general, if we take $f(\pi(s)) = i_{s}$ it is clear that $f$ is a $\pi$-partition if and only if $(i_1, i_2, \ldots, i_n)$ is a solution to \eqref{eq:is}, which gives \[\Omega_{\pi}(k) =\binom{k+n-1-\des(\pi)}{n}.\]

Before moving on, let us point out that in order to prove that the
formulas in this paper hold as polynomials in $x$ and $y$, it will
suffice to prove that they hold for all pairs of positive
integers. It is not hard to prove this fact, and we use it in each
of the proofs presented in this paper.

\subsection{Proofs for the case of the symmetric group}\label{sec:cycprf}
We will now prove Theorem \ref{thm:ges} using the theory of
$P$-partitions.

\begin{proof}[Proof of Theorem \ref{thm:ges}]
If we write out $\phi(xy)=\phi(x)\phi(y)$ using the
definition, we have
\begin{eqnarray*}
\sum_{\pi\in\mathfrak{S}_{n}}\binom{xy+n-1-\des(\pi)}{n}\pi & = &
\sum_{\sigma\in\mathfrak{S}_{n}}\binom{x+n-1-\des(\sigma)}{n}\sigma
\sum_{\tau\in\mathfrak{S}_{n}}\binom{y+n-1-\des(\tau)}{n}\tau \\
 & = & \sum_{\sigma,\tau \in \mathfrak{S}_n }\binom{x+n-1-\des(\sigma)}{n}\binom{y+n-1-\des(\tau)}{n}
\sigma\tau
\end{eqnarray*}
If we equate the coefficients of $\pi$ we have
\begin{equation}\label{eqst}
\binom{xy+n-1-\des(\pi)}{n} = \sum_{\sigma\tau = \pi}
\binom{x+n-1-\des(\sigma)}{n} \binom{y+n-1-\des(\tau)}{n}.
\end{equation}
Clearly, if formula \eqref{eqst} holds for all $\pi$, then formula
\eqref{eq1} is true. Let us interpret the left hand side of this
equation.

Let $x =k$, and $y=l$ be positive integers. Then the left hand
side of equation (\ref{eqst}) is just the order polynomial
$\Omega_{\pi}(kl)$. To compute this order polynomial we need to count the number of $\pi$-partitions $f:[n] \to X$, where $X$ is some totally ordered set with $kl$ elements. But instead of using $[kl]$ as our image set, we will use a different totally ordered set of the same cardinality. Let us count the $\pi$-partitions $f:[n]\to [l]\times[k]$. This is equal to the number of solutions to
\begin{equation}
(1,1) \leq (i_1,j_1) \leq (i_2,j_2) \leq \cdots \leq (i_n,j_n)
\leq (l,k)
\mbox{\quad and } (i_{s},j_s) < (i_{s+1},j_{s+1}) \mbox{ if } s \in \Des(\pi). \label{eqlex}
\end{equation}
Here we take the \emph{lexicographic
ordering} on pairs of integers. Specifically, $(i,j) < (i',j')$ if
$i < i'$ or else if $i = i'$ and $j < j'$.

To get the result we desire, we will sort the set of all solutions to \eqref{eqlex} into distinct cases indexed by subsets $I \subset [n-1]$. The sorting depends on $\pi$ and proceeds as follows. Let $F = ( (i_1, j_1), \ldots, (i_n, j_n) )$ be any solution to \eqref{eqlex}. For any $s=1,2,\ldots,n-1$, if $\pi(s) < \pi(s+1)$, then $(i_{s},j_{s}) \leq (i_{s+1},j_{s+1})$, which falls into one of
two mutually exclusive cases:
\begin{eqnarray}
i_{s} \leq i_{s+1} & \mbox{and} & j_{s}\leq j_{s+1}, \mbox{ or } \label{eqn:1}\\
i_{s} < i_{s+1} & \mbox{and} & j_{s} > j_{s+1}. \label{eqn:2}
\end{eqnarray}
If $\pi(s) > \pi(s+1)$, then $(i_{s},j_{s}) < (i_{s+1},j_{s+1})$,
which means either:
\begin{eqnarray}
i_{s} \leq i_{s+1} & \mbox{and} & j_{s} < j_{s+1}, \mbox{ or } \label{eqn:3}\\
i_{s} < i_{s+1} & \mbox{and} & j_{s} \geq j_{s+1},\label{eqn:4}
\end{eqnarray}
also mutually exclusive. Define $I_F = \{ s \in [n-1]\setminus \Des(\pi) \,\mid\, j_s > j_{s+1} \} \cup \{ s \in \Des(\pi) \,\mid\, j_s \geq j_{s+1} \}$. Then $I_F$ is the set of all $s$ such that either \eqref{eqn:2} or \eqref{eqn:4} holds for $F$. Notice that in both cases, $i_s < i_{s+1}$. Now for any $I \subset [n-1]$, let $S_I$ be the set of all solutions $F$ to \eqref{eqlex} satisfying $I_F = I$. We have split the solutions of \eqref{eqlex} into $2^{n-1}$ distinct cases indexed by all the different subsets $I$ of $[n-1]$.

Say $\pi = (2,1,3)$. Then we want to count the number of solutions to \[(1,1)\leq (i_{1},j_{1}) < (i_{2},j_{2}) \leq (i_{3},j_{3}) \leq (l,k),\] which splits into four distinct cases:
\begin{eqnarray*}
\emptyset\,\,: & i_{1} \leq i_{2} \leq i_{3} & j_{1} < j_{2} \leq j_{3} \\
\{\,1\,\}: & i_{1} < i_{2} \leq i_3 & j_1 \geq j_2 \leq j_3 \\
\{\,2\,\}: & i_1 \leq i_2 < i_3 & j_1 < j_2 > j_3 \\
\{\,1,2\,\}: & i_1 < i_2 < i_3 & j_1 \geq j_2 > j_3.
\end{eqnarray*}
We now want to count all the solutions contained in each of these cases and add them up. For a fixed subset $I$ we will use the theory of $P$-partitions to count the number of solutions for the set of inequalities first for the $j_{s}$'s and then for the $i_{s}$'s. Multiplying will give us the number of solutions in $S_I$; we do the same for the remaining subsets and sum to obtain the final result. For $I = \{\,1\,\}$ in the example above, we would count first the number of integer solutions to $j_1 \geq j_2 \leq j_3$, with $1\leq j_s \leq k$, and then we multiply this number by the number of solutions to $1 \leq i_1 < i_2 \leq i_3 \leq l$ to obtain the cardinality of $S_{\{1\}}$.  We will now carry out the computation in general.

For any particular $I\subset [n-1]$, form the poset $P_{I}$ of the elements
$1,2,\ldots,n$ by $\pi(s) <_{P_{I}} \pi(s+1)$ if $s \notin I$,
$\pi(s) >_{P_{I}} \pi(s+1)$ if $s\in I$. We form a ``zig-zag" poset of $n$ elements labeled consecutively by $\pi(1),
\pi(2),\ldots,\pi(n)$, with downward zigs corresponding to the
elements of $I$. For example, if $I =\{2,3\}$ for $n=5$, then $P_{I}$ has $\pi(1) <
\pi(2) > \pi(3) > \pi(4) < \pi(5)$.

\begin{figure} [h]
\centering
\includegraphics{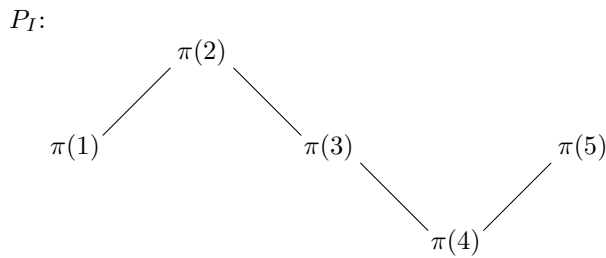}
\caption{The ``zig-zag" poset $P_{I}$ for $I = \{2,3\} \subset
[5]$.}
\end{figure}

For any solution in $S_I$, let $f: [n] \to [k]$ be defined by $f(\pi(s)) = j_{s}$ for $1\leq
s\leq n$. We will show that $f$ is a $P_{I}$-partition. If $\pi(s) <_{P_{I}} \pi(s+1)$ and $\pi(s) < \pi(s+1)$
in $\mathbb{Z}$, then \eqref{eqn:1} tells us that $f(\pi(s)) =
j_{s} \leq j_{s+1} = f(\pi(s+1))$. If $\pi(s) <_{P_{I}} \pi(s+1)$
and $\pi(s) > \pi(s+1)$ in $\mathbb{Z}$, then \eqref{eqn:3} tells
us that $f(\pi(s)) = j_{s} < j_{s+1} = f(\pi(s+1))$. If $\pi(s)
>_{P_{I}} \pi(s+1)$ and $\pi(s) < \pi(s+1)$ in $\mathbb{Z}$,
then \eqref{eqn:2} gives us that $f(\pi(s))= j_{s} > j_{s+1} =
f(\pi(s+1))$. If $\pi(s) >_{P_{I}} \pi(s+1)$ and $\pi(s) >
\pi(s+1)$ in $\mathbb{Z}$, then \eqref{eqn:4} gives us that
$f(\pi(s)) = j_{s} \geq j_{s+1} = f(\pi(s+1))$. In other words, we
have verified that $f$ is a $P_{I}$-partition. So for any particular solution in $S_I$, the $j_{s}$'s can be thought of as a $P_{I}$-partition. Conversely, any $P_{I}$-partition $f$ gives a solution in $S_I$ since if $j_s = f(\pi(s))$, then $(( i_1,j_1),\ldots,(i_n,j_n)) \in S_I$ if and only if $1 \leq i_1 \leq \cdots \leq i_n \leq l$ and $i_s < i_{s+1}$ for all $i \in I$. We can therefore turn our attention to counting $P_{I}$-partitions.

Let $\sigma \in \mathcal{L}(P_{I})$. Then for any
$\sigma$-partition $f$, we get a chain
\[1\leq f(\sigma(1)) \leq f(\sigma(2))\leq \cdots \leq f(\sigma(n)) \leq
k\] with $f(\sigma(s)) < f(\sigma(s+1))$ if $s \in \Des(\sigma)$.
The number of solutions to this set of inequalities is
\[\Omega_{\sigma}(k) = \binom{k+n-1-\des(\sigma)}{n}.\]

Recall by Observation \ref{ob1} that $\sigma^{-1}\pi(s) <
\sigma^{-1}\pi(s+1)$ if $\pi(s) <_{P_{I}} \pi(s+1)$, i.e., if $s \notin I$. If $\pi(s) >_{P_{I}} \pi(s+1)$ then
$\sigma^{-1}\pi(s) > \sigma^{-1}\pi(s+1)$ and $s \in I$. We get that $\Des(\sigma^{-1}\pi) = I$ if and only if $\sigma \in \mathcal{L}(P_I)$. Set $\tau =
\sigma^{-1}\pi$. The number of solutions to
\[1 \leq i_{1} \leq \cdots \leq i_{n} \leq l \mbox{ \quad and } i_{s} <
i_{s+1} \mbox{ if } s\in \Des(\tau)\] is given by
\[\Omega_{\tau}(l) = \binom{l+n-1-\des(\tau)}{n}.\] Now for a given $I$, the number of solutions in $S_I$ is
\[\sum_{\substack{\sigma\in\mathcal{L}(P_{I}) \\ \sigma\tau = \pi}}
\binom{k+n-1-\des(\sigma)}{n}\binom{l+n-1-\des(\tau)}{n}. \] Summing over all subsets $I \subset [n-1]$, we
can write the number of all solutions to \eqref{eqlex} as
\[\sum_{\sigma\tau=\pi}\binom{k+n-1-\des(\sigma)}{n}\binom{l+n-1-\des(\tau)}{n},
\] and so we have derived formula \eqref{eqst}.
\end{proof}

We now have a taste of how $P$-partitions can be used. We are
ready to go on and prove Theorem \ref{thmcyc}.

\begin{proof}[Proof of Theorem \ref{thmcyc}]
If we write out the definition for $\varphi(x)$ in the statement
of Theorem \ref{thmcyc}, multiply both sides by $n^2$, and equate coefficients, we have for any
$\pi \in \mathfrak{S}_{n}$, \[n\binom{xy+n-1-\cdes(\pi)}{n-1} =
\sum_{\sigma\tau = \pi} \binom{x+n-1-\cdes(\sigma)}{n-1}
\binom{y+n-1-\cdes(\tau)}{n-1}.\] For some $i$, we can write $\pi
= \nu\omega^{i}$ where $\omega$ is the $n$-cycle
$(\,1\,\,2\,\,\cdots \,\,n\,)$ and $\nu =
(n,\nu(2),\ldots,\nu(n))$. Observe that $\cdes(\pi) = \cdes(\nu) =
\des(\nu)$. Form the permutation $\widehat{\nu} \in
\mathfrak{S}_{n-1}$ by $\widehat{\nu}(s) = \nu(s+1)$, $s =
1,2,\ldots,n-1$. Then we can see that $\cdes(\pi) =
\des(\widehat{\nu}) + 1$. We have
\[ \binom{xy+n-1-\cdes(\pi)}{n-1}  =
\binom{xy+(n-1)-1-\des(\widehat{\nu})}{n-1}.
\]

Now we can apply equation \eqref{eqst} to give us
\begin{eqnarray}\label{eq:cyc}
\binom{xy+(n-1)-1-\des(\widehat{\nu})}{n-1}  \qquad \qquad \qquad \qquad \qquad \qquad \qquad \qquad \qquad & & \\
= \sum_{\sigma\tau =
\widehat{\nu}} \binom{x+(n-1)-1-\des(\sigma)}{n-1}\cdot\binom{y+(n-1)-1-\des(\tau)}{n-1}. & & \nonumber
\end{eqnarray}
For each pair of permutations $\sigma,\tau \in \mathfrak{S}_{n-1}$ such that $\sigma\tau
= \widehat{\nu}$, define the permutations $\widetilde{\sigma},
\widetilde{\tau} \in \mathfrak{S}_n$ as follows. For $s =
1,2,\ldots,n-1$, let $\widetilde{\sigma}(s) = \sigma(s)$ and
$\widetilde{\tau}(s+1) = \tau(s)$. Let $\widetilde{\sigma}(n) = n$
and $\widetilde{\tau}(1) = n$. Then by construction we have
$\widetilde{\sigma}\widetilde{\tau} = \nu$ and a quick observation
tells us that $\cdes(\widetilde{\sigma}) = \des(\sigma) + 1$ and
$\cdes(\widetilde{\tau}) = \des(\tau) + 1$. On the other hand, from any pair of permutations $\widetilde{\sigma}, \widetilde{\tau} \in \mathfrak{S}_n$ such that $\widetilde{\sigma}\widetilde{\tau} = \nu$, $\widetilde{\sigma}(n) = n$, we can construct a pair of permutations $\sigma, \tau \in \mathfrak{S}_{n-1}$ such that $\sigma\tau = \widehat{\nu}$ by reversing the process. Observe now that if $\widetilde{\sigma}(n) =n$ and $\widetilde{\sigma}\widetilde{\tau} = \nu$, then $\widetilde{\tau}(1) = n$. Therefore we have that
\eqref{eq:cyc} is equal to
\begin{eqnarray*}
& & \sum_{\substack{ \widetilde{\sigma}\widetilde{\tau} = \nu \\ \widetilde{\sigma}(n) = n }}
\binom{x+n-1-\cdes(\widetilde{\sigma})}{n-1}
\binom{y+n-1-\cdes(\widetilde{\tau})}{n-1}\\
& = & \sum_{\substack{ \widetilde{\sigma}(\widetilde{\tau}\omega^{i}) = \pi \\ \widetilde{\sigma}(n) = n } }
\binom{x+n-1-\cdes(\widetilde{\sigma})}{n-1}
\binom{y+n-1-\cdes(\widetilde{\tau}\omega^{i})}{n-1}\\
& = &
\sum_{\substack{ (\widetilde{\sigma}\omega^{n-j})(\omega^{j}\widetilde{\tau}\omega^{i}) =
\pi \\ \widetilde{\sigma}(n) = n } } \binom{x+n-1-\cdes(\widetilde{\sigma}\omega^{n-j})}{n-1}
\binom{y+n-1-\cdes(\omega^{j}\widetilde{\tau}\omega^{i})}{n-1}\\
& = & \sum_{\substack{ \sigma\tau = \pi \\ \sigma(j) = n} }
\binom{x+n-1-\cdes(\sigma)}{n-1} \binom{y+n-1-\cdes(\tau)}{n-1},
\end{eqnarray*}
where the last two formulas hold for any $j \in [n]$. Notice that the number of cyclic descents of $\tau =
\omega^{j}\widetilde{\tau}\omega^{i}$ is still the same as the
number of cyclic descents of $\widetilde{\tau}$. We take the sum over all $j = 1,\ldots, n$, yielding
\[n\binom{xy+n-1-\cdes(\pi)}{n-1} = \sum_{\sigma\tau = \pi}
\binom{x+n-1-\cdes(\sigma)}{n-1} \binom{y+n-1-\cdes(\tau)}{n-1}\]
as desired.
\end{proof}

\section{The augmented descent algebra}\label{sec:aug}
\subsection{Definitions and observations}\label{sec:def} Here we
present the definitions and basic results we will need to prove
the remaining theorems. Proofs of some of these basic facts are
identical to the proofs of analogous statements for ordinary
$P$-partitions and may be omitted. It bears mentioning that the
following definitions for type B posets and type B $P$-partitions, though taken from Chow \cite{Ch}, derive
from earlier work by Victor Reiner \cite{R}. In \cite{R2}, Reiner
extends this concept to any Coxeter group.
\begin{defn}
A \emph{$\mathfrak{B}_{n}$ poset} is a poset $P$ whose elements
are $0, \pm 1, \pm 2, \ldots, \pm n$ such that if $i <_{P} j$ then
$-j <_{P} -i$. Note that if we are given a poset with $n+1$
elements labeled by $0, a_{1},\ldots, a_{n}$ where $a_{i} = i$ or
$-i$, then we can extend it to a $\mathfrak{B}_{n}$ poset of
$2n+1$ elements.
\end{defn}

\begin{figure} [h]
\centering
\includegraphics[scale=.8]{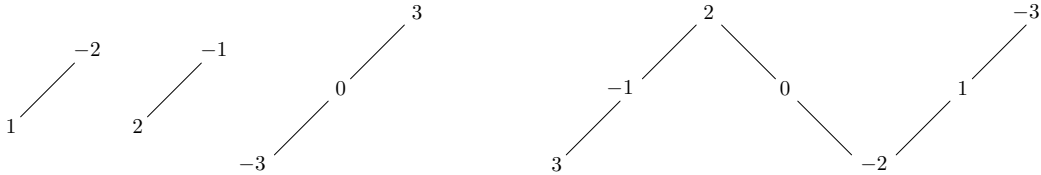}
\caption{Two $\mathfrak{B}_{3}$ posets.}
\end{figure}

Let $X = \{x_{0},x_{1},x_2, \ldots\}$ be a countable, totally ordered set with total order \[ x_0 < x_1 < x_2 < \cdots .\] Then define $\pm X$ to be the set $\{\ldots,-x_{1},x_{0},x_{1},\ldots\}$ with total order \[ \cdots < -x_2 < -x_1 < x_0 < x_1 < x_2 < \cdots .\]

\begin{defn}\label{def:typebppart} For
any $\mathfrak{B}_{n}$ poset $P$, a \emph{$P$-partition of type B} is a function $f:\pm [n] \to \pm X$ such that:
\begin{itemize}
\item $f(i) \leq f(j)$ if $i <_{P} j$

\item $f(i) < f(j)$ if $i <_{P} j$ and $i > j$ in $\mathbb{Z}$

\item $f(-i) = -f(i)$
\end{itemize}
\end{defn}

Note that type B $P$-partitions differ from ordinary $P$-partitions only in the addition of the property $f(-i) = -f(i)$. Let $\mathcal{A}(P)$ denote the set of all type B $P$-partitions. We usually think of $X$ as a subset of the nonnegative integers, and when $X$ has finite cardinality $k+1$, then the \emph{type B order polynomial}, denoted $\Omega_{P}(k)$, is the number of $P$-partitions $f:\pm[n] \to \pm X$.

We can think of any signed permutation $\pi \in \mathfrak{B}_{n}$ as a $\mathfrak{B}_{n}$ poset with the total order $\pi(s) <_{\pi} \pi(s+1)$, $0\leq s
\leq n-1$. For example, the signed permutation $(-2,1)$ has $-1
<_{\pi} 2 <_{\pi} 0 <_{\pi} -2 <_{\pi} 1$ as a poset. Note that
$\mathcal{A}(\pi)$ is the set of all functions $f: \pm[n] \to \pm
X$ such that for $0 \leq s \leq n$, $f(-s) = -f(s)$ and
\[ x_0 = f(\pi(0)) \leq f(\pi(1)) \leq f(\pi(2))\leq \cdots \leq
f(\pi(n)),\]
where $f(\pi(s))< f(\pi(s+1))$ whenever $\pi(s) > \pi(s+1)$, $s = 0,1,\ldots,n-1$.
For example, the type B $\pi$-partitions where $\pi = (-2,1)$ are all maps $f$
such that $x_0 < f(-2) \leq f(1)$. As before,
\begin{obs}\label{ob2} For a type B poset $P$, a signed permutation $\pi$ is in $\mathcal{L}(P)$ if and only if $i<_{P}
j$ implies $\pi^{-1}(i) < \pi^{-1}(j)$.
\end{obs}
We have a fundamental theorem for $P$-partitions of type B.
\begin{thm}[FTPPB]
The set of all type B $P$-partitions of a $\mathfrak{B}_{n}$
poset $P$ is the disjoint union of the set of $\pi$-partitions of
all linear extensions $\pi$ of $P$: \[\mathcal{A}(P) =
\coprod_{\pi \in \mathcal{L}(P)} \mathcal{A}(\pi).\]
\end{thm}
\begin{cor}
\[\Omega_{P}(k) = \sum_{\pi \in \mathcal{L}(P)} \Omega_{\pi}(k).\]
\end{cor}

\begin{figure} [h]
\centering
\includegraphics{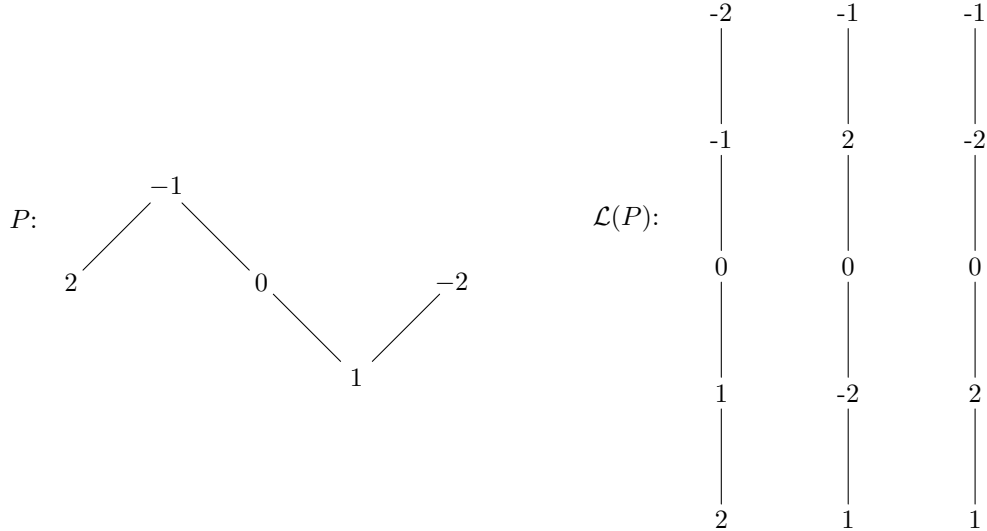}
\caption{Linear extensions of a $\mathfrak{B}_{2}$ poset $P$.}
\end{figure}

Similarly to the type A case, it is easy to compute the order polynomial $\Omega_{\pi}(k)$ for any permutation $\pi \in \mathfrak{B}_n$. Any $\pi$-partition $f: \pm[n] \to \pm[k]$ is determined by where we map $\pi(1), \pi(2),\ldots,\pi(n)$ since $f(-i) = -f(i)$ (and so $f(0) = 0$). To count them we can take $f(\pi(s)) = i_s$, and look at the number of integer solutions to the set
of inequalities \[ 0 \leq i_1 \leq i_2 \leq \cdots \leq i_n \leq
k,  \mbox{ \quad and  } i_s < i_{s+1} \mbox{ if } s \in \Des(\pi).\] This number is
the same as the number of solutions to \[ 1\leq i_1 \leq i_2 \leq
\cdots \leq i_n \leq k+1-\des(\pi),\] which we know to be $\big(\!\binom{k+1 -
\des(\pi)}{n}\!\big)$. We have \[\Omega_{\pi}(k) = \binom{k+n-\des(\pi)}{n}.\]

Now we have the tools to prove Theorem \ref{thm:cho}.
\begin{proof}[Proof of Theorem \ref{thm:cho}]
We will omit the details, since they are essentially the same as for Theorem \ref{thm:ges}. The main difference is that we want to count $\pi$-partitions $f: \pm[n] \to \pm[l]\times\pm[k]$. We notice that because of the property $f(-s) = -f(s)$, this is just like counting all $f: [n] \to \{0,1,\ldots,l\}\times \{-k,\ldots,-1,0,1,\ldots, k\}$ where for $s=1,2,\ldots,n$,  $f(\pi(s)) = (i_s ,j_s )$ with $(0,0) \leq (i_s,j_s) \leq (l,k)$ in the lexicographic order.

\begin{figure} [h]
\centering
\includegraphics[scale=.8]{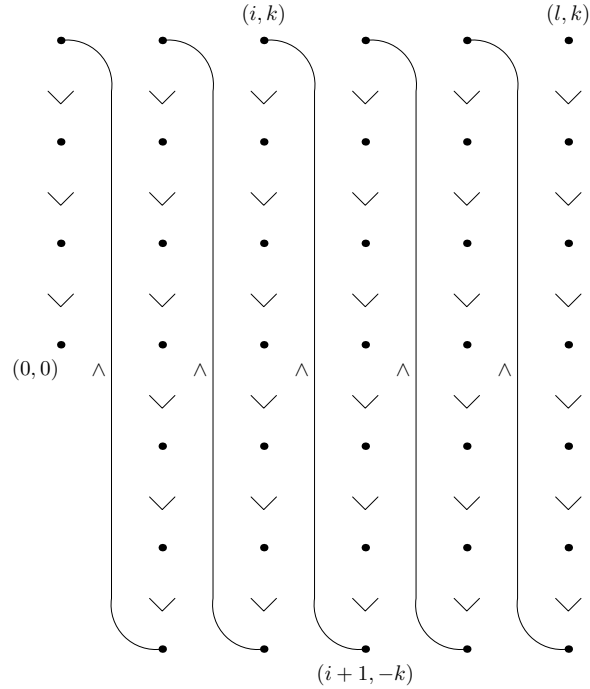}
\caption{The lexicographic order on $\{0,1,\ldots,l\}\times \{-k,\ldots,-1,0,1,\ldots, k\}$ with $(0,0) \leq (i_s,j_s) \leq (l,k)$.}
\end{figure}

This choice of image set $X$ has $2kl + k + l + 1$ elements, and so for each $\pi$ we can count all these maps with $\Omega_{\pi}(2kl+k+l) = \binom{2kl+k+l - \des(\pi)}{n}$. We use similar arguments to those of Theorem \ref{thm:ges} for splitting the lexicographic solutions to \[(0,0) \leq (i_1,j_1) \leq \cdots \leq (i_n,j_n) \leq (l,k) \mbox{ \quad and } (i_s,j_s) < (i_{s+1},j_{s+1}) \mbox{ if } s\in \Des(\pi).\] Once we have properly grouped the set of solutions it is not much more work to obtain the crucial formula: \[\Omega_{\pi}(2kl+k+l) = \sum_{\sigma\tau = \pi} \Omega_{\sigma}(k)\Omega_{\tau}(l).\]
\end{proof}

We now give the definition of an augmented $P$-partition and basic tools related to their study. Let $X = \{x_{0},x_{1},\ldots,x_{\infty}\}$ be a countable,
totally ordered set with a maximal element $x_{\infty}$. The total
ordering on $X$ is given by \[x_{0} < x_{1} < x_2 < \cdots < x_{\infty}.\] Define $\pm X$ to be
$\{-x_{\infty},\ldots,-x_{1},x_{0},x_{1},\ldots,x_{\infty}\}$ with the total order \[-x_{\infty} < \cdots < -x_1 < x_0 < x_1 < \cdots < x_{\infty}.\]

\begin{defn}\label{def:augp}
For any $\mathfrak{B}_{n}$ poset $P$, an \emph{augmented
$P$-partition} is a function $f:\pm [n] \to \pm X$ such that:
\begin{itemize}
\item $f(i) \leq f(j)$ if $i <_{P} j$

\item $f(i) < f(j)$ if $i <_{P} j$ and $i > j$ in $\mathbb{Z}$

\item $f(-i) = -f(i)$

\item if $0 < i$ in $\mathbb{Z}$, then $f(i) <
x_{\infty}$.
\end{itemize}
\end{defn}
Note that augmented $P$-partitions differ from $P$-partitions of
type B only in the addition of maximal and minimal elements of the image set $\pm X$ and in the last criterion. Let
$\mathcal{A}^{(a)}(P)$ denote the set of all augmented $P$-partitions.
When $X$ has finite cardinality $k+1$ (and so $\pm X$ has cardinality $2k+1$), then the \emph{augmented order
polynomial}, denoted $\Omega_{P}^{(a)}(k)$, is the number of
augmented $P$-partitions.

For any signed permutation $\pi \in \mathfrak{B}_n$, note that $\mathcal{A}^{(a)}(\pi)$ is the set of all functions $f: \pm[n] \to \pm X$ such that for $0 \leq s \leq n$, $f(-s) = -f(s)$ and
\[x_0 = f(\pi(0)) \leq f(\pi(1)) \leq f(\pi(2))\leq \cdots \leq
f(\pi(n)) \leq x_{\infty},\] and $f(\pi(s))< f(\pi(s+1))$ whenever $\pi(s) > \pi(s+1)$. In addition, we
have $f(\pi(n)) < x_{\infty}$ whenever $\pi(n) > 0$. The set of all
augmented $\pi$-partitions where $\pi = (-2,1)$ is all maps $f$
such that $x_0 < f(-2) \leq f(1) < x_{\infty}$.

We have a fundamental theorem for augmented $P$-partitions.
\begin{thm}[FTAPP]
The set of all augmented $P$-partitions of a $\mathfrak{B}_{n}$
poset $P$ is the disjoint union of the set of $\pi$-partitions of
all linear extensions $\pi$ of $P$: \[\mathcal{A}^{(a)}(P) =
\coprod_{\pi \in \mathcal{L}(P)} \mathcal{A}^{(a)}(\pi).\]
\end{thm}
\begin{cor}
\[\Omega^{(a)}_{P}(k) = \sum_{\pi \in \mathcal{L}(P)} \Omega^{(a)}_{\pi}(k).\]
\end{cor}
As in our previous cases, it is fairly easy to compute the augmented order polynomial for totally ordered chains. The number of augmented $\pi$-permutations $f: \pm[n]
\to \pm[k]$ is equal to the number of integer solutions to the set
of inequalities \[ 0 \leq i_1 \leq i_2 \leq \cdots \leq i_n \leq
k = i_{n+1},  \mbox{ \quad and  } i_s < i_{s+1} \mbox{ if } s \in \aDes(\pi).\] Therefore, just as with other types of order polynomials, we can express the augmented order polynomial as a binomial coefficient,
\[\Omega^{(a)}_{\pi}(k) = \binom{ k + n - \ades(\pi)}{n}.\]

To give one example of the usefulness of augmented $P$-partitions,
we conclude this subsection with proof of Proposition \ref{prp:augeul}, claiming $A_{n}^{(a)}(t) = 2^{n}A_{n}(t)$.

\begin{proof}[Proof of Proposition \ref{prp:augeul}]
From the general theory of $P$-partitions in Stanley's book
\cite{St}, we have
\[\sum_{k\geq 0}\Omega_{P}(k)t^{k} = \frac{\sum_{\pi\in
\mathcal{L}(P)} t^{\des(\pi)+1}}{(1-t)^{|P|+1}}. \] Let $P$ be an
antichain---that is, a poset with no relations---of $n$ elements.
Then $\Omega_{P}(k) = k^n$ since each of the $n$ elements of $P$
is free to be mapped to any of $k$ places. Furthermore,
$\mathcal{L}(P) = \mathfrak{S}_{n}$, so we get the following
equation,
\[\sum_{k\geq 0}k^{n}t^{k} = \frac{A_{n}(t)}{(1-t)^{n+1}}.\]

Now let $P$ be the poset given by an antichain of $2n+1$ elements
labeled $0,\pm 1, \pm 2,\ldots, \pm n$. The number of augmented
$P$-partitions $f: \pm[n] \to \pm[k]$ is determined by the choices
for $f(1),f(2),\ldots,f(n)$, which can take any of the values in
$\{-k,-k+1,\ldots,k-2,k-1\}$. Therefore $\Omega^{(a)}_{P}(k) =
(2k)^{n}$. For $\mathfrak{B}_n$ posets $P$, it is not difficult to
show that we have the identity \[\sum_{k\geq
0}\Omega^{(a)}_{P}(k)t^{k} = \frac{\sum_{\pi\in \mathcal{L}(P)}
t^{\ades(\pi)}}{(1-t)^{n+1}}.
\] For our antichain we have $\mathcal{L}(P) = \mathfrak{B}_{n}$,
and therefore
\[\frac{A^{(a)}_{n}(t)}{(1-t)^{n+1}} = \sum_{k\geq 0}(2k)^{n}t^{k}
= 2^{n}\sum_{k\geq 0}k^{n}t^{k} =
\frac{2^{n}A_{n}(t)}{(1-t)^{n+1}},\] so the proposition is proved.
\end{proof}

\subsection{Proofs for the case of the hyperoctahedral group}\label{sec:prf}
The following proofs will follow the same basic structure as the
proof of Theorem \ref{thm:ges}, but with some important changes in
detail. In both cases we will rely on a slightly different total
ordering on the integer points $(i,j)$, where $i$ and $j$ are
bounded both above and below. Let us now define the
\emph{augmented lexicographic order}.

\begin{figure} [h]
\centering
\includegraphics[scale=.7]{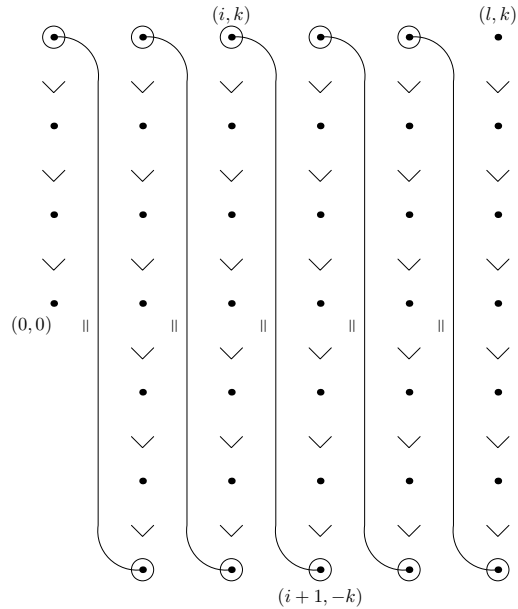}
\caption{The augmented lexicographic order.}
\end{figure}

Consider all points $(i,j)$ with $0 \leq i \leq l$, $-k \leq j
\leq k$. We have $(i,j) < (i',j')$ if $i< i'$ or else if $i = i'$
and $j < j'$ as before, except in the important special case that
follows. We now say $(i,j) = (i',j')$ in one of two situations.
Either
\begin{eqnarray*}
i = i' & \textrm{ and } & j = j'
\end{eqnarray*}
or
\begin{eqnarray*}
 i+1 = i' & \textrm{ and } & j = k = -j'.
\end{eqnarray*}
If we have $0\leq i \leq l$, $-2 \leq j \leq 2$, then in augmented
lexicographic order, the first few points $(0,0) \leq (i,j) \leq
(l,2)$ are:

\begin{eqnarray*}
& & (0,0) < (0,1) < (0,2) \\
 & & \qquad \qquad \qquad =  (1,-2) < (1,-1) < (1,0) < (1,1) < (1,2) \\
 & & \qquad \qquad \qquad \qquad \qquad \qquad = (2,-2) < (2,-1) < (2,0) < \cdots
\end{eqnarray*}

To be more precise, what we have done is to form equivalence
classes of points and to introduce a total order on these
equivalence classes. If $j\neq \pm k$, then the class represented
by $(i,j)$ is just the point itself. Otherwise, the classes
consist of the two points $(i,k)$ and $(i+1,-k)$. When we write
$(i,j) = (i',j')$, what we mean is that the two points are in the
same equivalence class. In the proofs that follow, it will be
important to remember the original points as well as the
equivalence classes to which they belong. This special ordering will be very apparent in deriving the $q$-analogs of Theorem \ref{conj1} and Theorem \ref{thm:hypermod}.

We will now prove
Theorem \ref{conj1}.

\begin{proof}[Proof of Theorem \ref{conj1}]
As before, we equate coefficients and prove that a simpler
formula,
\begin{eqnarray}\label{eq:c1}
\binom{2kl+n-\ades(\pi)}{n} \!\!\! & = \!\!\! & \sum_{\sigma\tau =
\pi}\binom{k+n-\ades(\sigma)}{n}\binom{l+n-\ades(\tau)}{n},
\end{eqnarray}
holds for any $\pi \in \mathfrak{B}_{n}$.

We recognize the left-hand side of equation \eqref{eq:c1} as
$\Omega_{\pi}^{(a)}(2kl)$, so we want to count augmented
$P$-partitions $f:\pm[n]\to \pm X$, where $X$ is a totally ordered
set of order $2kl+1$. We interpret this as the number of
solutions, in the augmented lexicographic ordering, to
\begin{equation}\label{eq:augin}
(0,0) \leq (i_{1},j_1) \leq (i_2,j_2) \leq \cdots \leq (i_n,j_n)
\leq (l,0),
\end{equation}
where we have
\begin{itemize}
\item $0\leq i_s \leq l$,
\item $-k<j_{s}\leq k$ if $\pi(s) <0$,
\item $-k\leq j_{s} < k$ if $\pi(s) > 0$, and
\item $(i_s,j_s) < (i_{s+1},j_{s+1})$ if $s \in \aDes(\pi)$.
\end{itemize}
Let us clarify. There are $2kl+l+1$ points $(i,j)$ with $0\leq
i \leq l$ and $-k \leq j \leq k$, not including the points
$(0,j)$ with $j<0$, or the points $(l,j)$ with $j > 0$. Under
the augmented lexicographic ordering, $l$ of these points are
identified: points of the form $(i,k) = (i+1,-k)$, for $i = 0,
1,\ldots, l-1$. Any particular $(i_{s},j_s)$ can only occupy one of
$(i,k)$ or $(i+1,-k)$, but not both. So there are truly
$2kl+1$ distinct classes in which the $n$ points can fall.
This confirms our interpretation of the order polynomial.

Now as before, we will split the solutions to the inequalities into distinct cases.
Let $\pi(0) = \pi(n+1) = 0$, $i_0 = j_0 = 0$, $i_{n+1} = l$, and
$j_{n+1} = 0$. Let $F = ( (i_1, j_1), \ldots, (i_n, j_n) )$ be any solution to \eqref{eq:augin}. If $\pi(s) < \pi(s+1)$, then $(i_{s},j_{s}) \leq
(i_{s+1},j_{s+1})$, which falls into one of two mutually exclusive
cases:
\begin{eqnarray}
i_{s} \leq i_{s+1} & \mbox{and} & j_{s}\leq j_{s+1}, \mbox{ or } \label{eqn:a}\\
i_{s} < i_{s+1} & \mbox{and} & j_{s} > j_{s+1}. \label{eqn:b}
\end{eqnarray}
If $\pi(s) > \pi(s+1)$, then $(i_{s},j_{s}) < (i_{s+1},j_{s+1})$,
which we split as:
\begin{eqnarray}
i_{s} \leq i_{s+1} & \mbox{and} & j_{s} < j_{s+1}, \mbox{ or } \label{eqn:c}\\
i_{s} < i_{s+1} & \mbox{and} & j_{s} \geq j_{s+1},\label{eqn:d}
\end{eqnarray}
also mutually exclusive. Define $I_F = \{ s \in \{0,1,\ldots,n\} \setminus \aDes(\pi) \,\mid\, j_s > j_{s+1} \} \cup \{ s \in \aDes(\pi) \,\mid\, j_s \geq j_{s+1} \}$. Then $I_F$ is the set of all $s$ such that either \eqref{eqn:b} or \eqref{eqn:d} holds for $F$. Now for any $I \subset \{0,1,\ldots,n\}$, let $S_I$ be the set of all solutions $F$ to \eqref{eq:augin} satisfying $I_F = I$. We have split the solutions of \eqref{eq:augin} into $2^{n+1}$ distinct cases indexed by all the different subsets $I$ of $\{0,1,\ldots,n\}$.

However, $S_{\emptyset}$ is empty, since \[0\leq i_{1} \leq \cdots \leq i_{n} \leq
l\] yields \[0\leq j_{1} \leq \cdots \leq j_{n} \leq 0 \mbox{ \quad with } j_{s} < j_{s+1} \mbox{ if } s \in \aDes(\pi).\] As discussed
before, the augmented descent set of a signed permutation is never
empty, so we would get $0 < 0$, a contradiction. At the other
extreme, the set $S_{\{0,1,\ldots,n\}}$ has no solutions either. Here we get \[0 < i_{1}
< \cdots < i_{n} < l\] and consequently \[0\geq j_{1} \geq
\cdots \geq j_{n} \geq 0 \mbox{ \quad with } j_{s} > j_{s+1} \mbox{ if } s
\notin \aDes(\pi).\] But $\aDes(\pi)$ cannot equal
$\{0,1,\ldots,n\}$, so we get the contradiction $0> 0$.

Now let $I$ be any nonempty, proper subset of $\{0,1,\ldots,n\}$.
Form the poset $P_{I}$ by $\pi(s) >_{P_{I}} \pi(s+1)$ if $s \in
I$, $\pi(s) <_{P_{I}} \pi(s+1)$ otherwise. The poset $P_{I}$ looks
like a zig-zag, labeled consecutively by
$0=\pi(0),\pi(1),\pi(2),\ldots,\pi(n),0=\pi(n+1)$ with downward
zigs corresponding to the elements of $I$. Because $I$ is neither
empty nor full, we never have $0<_{P_{I}} 0$, so $P_{I}$ is a
well-defined, nontrivial type B poset.

For a given $F \in S_I$, let $f: \pm[n] \to \pm [k]$ be defined by $f(\pi(s)) = j_{s}$ and $f(-s) = -f(s)$ for $s = 0,1,\ldots,n$. We will show that $f$ is an augmented $P_{I}$ partition. If $\pi(s) <_{P_{I}} \pi(s+1)$ and $\pi(s) < \pi(s+1)$ in $\mathbb{Z}$, then \eqref{eqn:a} tells us that $f(\pi(s)) = j_{s} \leq j_{s+1} =
f(\pi(s+1))$. If $\pi(s) <_{P_{I}} \pi(s+1)$ and $\pi(s) >
\pi(s+1)$ in $\mathbb{Z}$, then \eqref{eqn:c} tells us that
$f(\pi(s)) = j_{s} < j_{s+1} = f(\pi(s+1))$. If $\pi(s)
>_{P_{I}} \pi(s+1)$ and $\pi(s) < \pi(s+1)$ in $\mathbb{Z}$,
then \eqref{eqn:b} gives us that $f(\pi(s))= j_{s} > j_{s+1} =
f(\pi(s+1))$. If $\pi(s) >_{P_{I}} \pi(s+1)$ and $\pi(s)
> \pi(s+1)$ in $\mathbb{Z}$, then \eqref{eqn:d} gives us that $f(\pi(s)) =
j_{s} \geq j_{s+1} = f(\pi(s+1))$. Since we required that
$-k<j_{s}\leq k$ if $\pi(s) < 0$ and $-k\leq j_{s} < k$ if
$\pi(s) > 0$, we have that for any particular solution in $S_I$, the $j_{s}$'s can be thought of as an augmented $P_{I}$-partition. Conversely, any augmented $P_{I}$-partition $f$ gives a solution in $S_I$ since if $j_s = f(\pi(s))$, then $(( i_1,j_1),\ldots,(i_n,j_n)) \in S_I$ if and only if $0 \leq i_1 \leq \cdots \leq i_n \leq l$ and $i_s < i_{s+1}$ for all $i \in I$. We can therefore turn our attention to counting augmented $P_{I}$-partitions.

Let $\sigma \in \mathcal{L}(P_{I})$. Then we get for any
$\sigma$-partition $f$,
\[0 \leq f(\sigma(1)) \leq f(\sigma(2))\leq \cdots \leq f(\sigma(n)) \leq
k,\] and $f(\sigma(s)) < f(\sigma(s+1))$ whenever $s \in
\aDes(\sigma)$, where we take $f(\sigma(n+1)) = k$. The number
of solutions to this set of inequalities is
\[\Omega^{(a)}_{\sigma}(k) = \binom{k+n-\ades(\sigma)}{n}.\]

Recall by Observation \ref{ob2} that $\sigma^{-1}\pi(s) <
\sigma^{-1}\pi(s+1)$ if $\pi(s) <_{P_{I}} \pi(s+1)$, i.e., if
$s\notin I$. If $\pi(s) >_{P_{I}} \pi(s+1)$ then
$\sigma^{-1}\pi(s) > \sigma^{-1}\pi(s+1)$ and $s\in I$.
We get that $\aDes(\sigma^{-1}\pi) = I$ if and only if $\sigma \in \mathcal{L}(P_I)$. Set $\tau =
\sigma^{-1}\pi$. The number of solutions to
\[0 \leq i_{1} \leq \cdots \leq i_{n} \leq l \mbox{ \quad and } i_{s} <
i_{s+1} \mbox{ if } s\in \aDes(\tau)\] is given by
\[\Omega_{\tau}(l) = \binom{l+n-\ades(\tau)}{n}.\] Now for a given $I$, the number of solutions to \eqref{eq:augin} is
\[\sum_{\substack{\sigma\in\mathcal{L}(P_{I}) \\ \sigma\tau = \pi}}
\binom{k+n-\ades(\sigma)}{n}\binom{l+n-\ades(\tau)}{n}. \] Summing over all subsets $I \subset \{0,1,\ldots,n\}$, we can write the number of all solutions to \eqref{eq:augin} as
\[\sum_{\sigma\tau=\pi}\binom{k+n-\ades(\sigma)}{n}\binom{l+n-\ades(\tau)}{n},
\] and so the theorem is proved.

\end{proof}

The proof of Theorem \ref{thm:hypermod} proceeds nearly identically, so
we will omit unimportant details in the proof below.

\begin{proof}[Proof of Theorem \ref{thm:hypermod}]
We equate coefficients and prove that
\begin{eqnarray}\label{eq:c2}
\binom{2kl+k+n-\ades(\pi)}{n} \!\!\!&=\!\!\!& \sum_{\sigma\tau =
\pi}\binom{k+n-\ades(\sigma)}{n}\binom{l+n-\des(\tau)}{n},
\end{eqnarray}
holds for any $\pi \in \mathfrak{B}_{n}$.

We recognize the left-hand side of equation \eqref{eq:c2} as
$\Omega_{\pi}^{(a)}(2kl+k)$, so we want to count augmented
$P$-partitions $f:\pm[n]\to \pm X$, where $X$ is a totally ordered
set of order $2kl+k+1$. We interpret this as the number of
solutions, in the augmented lexicographic ordering, to
\begin{equation}\label{eq:augin2}
(0,0) \leq (i_{1},j_1) \leq (i_2,j_2) \leq \cdots \leq (i_n,j_n)
\leq (l,k),
\end{equation}
where we have
\begin{itemize}
\item $0\leq i_s \leq l$,
\item $-k<j_{s}\leq k$ if $\pi(s) <0$,
\item $-k\leq j_{s} < k$ if $\pi(s) > 0$, and
\item $(i_s,j_s) < (i_{s+1},j_{s+1})$ if $s \in \aDes(\pi)$.
\end{itemize}
With these restrictions, we split the solutions to \eqref{eq:augin2} by our prior rules. Let $F = ( (i_1, j_1), \ldots , (i_n, j_n) )$ be any particular solution. If $\pi(s) < \pi(s+1)$, then
$(i_{s},j_{s}) \leq (i_{s+1},j_{s+1})$, which falls into one of two mutually exclusive cases:
\begin{eqnarray*}
i_{s} \leq i_{s+1} & \mbox{and} & j_{s}\leq j_{s+1}, \mbox{ or }\\
i_{s} < i_{s+1} & \mbox{and} & j_{s} > j_{s+1}.
\end{eqnarray*}
If $\pi(s) > \pi(s+1)$, then $(i_{s},j_{s}) < (i_{s+1},j_{s+1})$,
giving:
\begin{eqnarray*}
i_{s} \leq i_{s+1} & \mbox{and} & j_{s} < j_{s+1}, \mbox{ or } \\
i_{s} < i_{s+1} & \mbox{and} & j_{s} \geq j_{s+1},
\end{eqnarray*}
also mutually exclusive. With $(i_n,j_n)$, there is only one case, depending on $\pi$. If $\pi(n) > 0$, then $(i_{n},j_{n})
< (l,k)$ and $i_{n}\leq l$ and $-k\leq j_{n}<k$. Similarly,
if $\pi(n) < 0$, then $(i_{n},j_{n})\leq (l,k)$ and we have
$i_{n} \leq l$ and $-k<j_{n}\leq k$. Define $I_F$ and $S_I$ as before. We get $2^n$ mutually
exclusive sets $S_I$ indexed by subsets $I \subset \{0,1,\ldots, n-1\}$ (these subsets will correspond to ordinary descent sets).

Now for any $I \subset \{0,1,\ldots,n-1\}$, define the $\mathfrak{B}_{n}$ poset $P_{I}$ to be the poset given by $\pi(s)
>_{P_{I}} \pi(s+1)$ if $s\in I$, and $\pi(s) <_{P_{I}}
\pi(s+1)$ if $s\notin I$, for $s = 0,1,\ldots,n-1$. We form a
zig-zag poset labeled consecutively by $\pi(0) = 0, \pi(1),
\pi(2),\ldots,\pi(n)$.

For a given solution $F \in S_I$, let $f: \pm[n] \to \pm[k]$ be defined by $f(\pi(s)) = j_{s}$ for
$0\leq s\leq n$, with $f(-s) = -f(s)$. It is not too difficult to
check that $f$ is an augmented $P_{I}$-partition, and that any augmented $P_I$-partition corresponds to a solution in $S_I$. Let $\sigma \in
\mathcal{L}(P_{I})$.  Then for any $\sigma$-partition $f$ we get
\[f(\sigma(0))=0 \leq f(\sigma(1))\leq \cdots \leq f(\sigma(n)) \leq
k,\] with $f(\sigma(s)) < f(\sigma(s+1))$ whenever $s \in
\aDes(\sigma)$. The number of solutions to this set of
inequalities is
\[\Omega^{(a)}_{\sigma}(k) = \binom{k+n-\ades(\sigma)}{n}.\]

We see that for $s = 0,1,\ldots, n-1$, $\sigma^{-1}\pi(s) <
\sigma^{-1}\pi(s+1)$ if $\pi(s) <_{P_{I}} \pi(s+1)$, i.e., if
$s\notin I$. Also, if $\pi(s) >_{P_{I}} \pi(s+1)$ then
$\sigma^{-1}\pi(s) > \sigma^{-1}\pi(s+1)$ and $s\in I$.
This time we get that $\Des(\sigma^{-1}\pi) = I$, an ordinary
descent set, if and only if $\sigma \in \mathcal{L}_{P_I}$. Set $\tau = \sigma^{-1}\pi$. The number of solutions
to \[0\leq i_{1} \leq \cdots \leq i_{n} \leq l \mbox{ \quad and } i_{s} < i_{s+1} \mbox{ if } s\in \Des(\tau)\] is given by
\[\Omega_{\tau}(l) = \binom{l+n-\des(\tau)}{n}.\] We take the sum over all subsets
$I$ to show the number of solutions to \eqref{eq:c2} is \[\sum_{\sigma\tau=\pi}\binom{k+n-\ades(\sigma)}{n}\binom{l+n-\des(\tau)}{n}, \] and the theorem is proved.
\end{proof}

\section{Some $q$-analogs}

In this section we give some $q$-analogs of our previous theorems. First we define the $q$-variant of the (type A) order polynomial, or $q$-order polynomial.\footnote{Properly speaking, this $q$-analog of the order polynomial is not a polynomial in $k$. However, we will refer to it as the ``$q$-order polynomial," even if it might be more appropriate to call it the ``$q$-analog of the order polynomial." } Let $P$ be a poset with $n$ elements, and take $X = \{0,1,\ldots,k-1\}$. Then, \[\Omega_{P}(q;k) = \sum_{f \in \mathcal{A}(P) } \left(\prod_{i =1}^{n} q^{f(i)}\right).\]
Let $n_{q}! = (1+q)(1+q+q^2)\cdots(1+q+\cdots+ q^{n-1})$ and define the \emph{$q$-binomial coefficient} $\binom{a}{b}_{q}$ in the natural way: \[\binom{a}{b}_q = \frac{a_q !}{b_q !(a-b)_q !}\] Another way to interpret the $q$-binomial coefficient is as the coefficient of $x^{b}y^{a-b}$ in $(x+y)^{a}$ where $x$ and $y$ ``$q$-commute" via the relation $yx = qxy$. This interpretation is good for some purposes, but we will use a different point of view. We will define the $q$-multi-choose function $\big(\!\binom{a}{b}\!\big)_q = \binom{a+b-1}{b}_q$ as the following: \[\left(\!\!\binom{a}{b}\!\!\right)_q = \sum_{0\leq i_{1} \leq \cdots \leq i_{b} \leq a-1} \left( \prod_{s =1}^{n} q^{i_s}\right) = \Omega_{id} (q;a),\] where $id$ is the totally ordered chain $1 <_P 2 <_P \cdots <_P b$, the identity permutation in $\mathfrak{S}_b$. In this case, the $q$-order polynomial counts certain integer partitions (so there is at least some overlap between $P$-partitions and the more familiar integer partitions). Specifically, the coefficient of $q^r$ is the number of integer partitions of $r$ with at most $b$ parts, each of size at most $a-1$.\footnote{It appears Fulman may have been aware of a $q$-analog of his shuffling results, though he never states $q$-analogs of his main theorems. See \cite{F2}, Definitions 3 and 5.}

Now we will obtain formulas for $q$-order polynomials for all permutations. For any permutation $\pi \in \mathfrak{S}_n$, the $q$-order polynomial may be expressed as \[\Omega_{\pi}(q;k) = \sum_{\substack{0\leq i_1 \leq \cdots \leq i_n \leq k-1 \\ s\in \Des(\pi) \Rightarrow i_s < i_{s+1}}} \left(\prod_{s=1}^{n} q^{i_s}\right).\] When we computed the ordinary order polynomial in section \ref{sec:def} we only cared about the \emph{number} of solutions, rather than the \emph{set} of solutions to the inequalities
\begin{equation}\label{eq:qineq}
0\leq i_1 \leq \cdots \leq i_n \leq k-1 \mbox{ \quad and } i_{s} < i_{s+1} \mbox{ if } s\in \Des(\pi).
\end{equation}
Since we only cared how many solutions there were and not what the solutions were, we could count solutions to a system where all the inequalities were weak. We will still follow the same basic procedure, but as we manipulate our system of inequalities we need to keep track of how we modify the set of solutions. The $q$-order polynomial will be seen to be simply a power of $q$ (depending on $\pi$) times a $q$-binomial coefficient.

We can form a new system of inequalities that has the same number of solutions as \eqref{eq:qineq}, but in which every inequality is weak: \[0\leq i'_1 \leq \cdots \leq i'_n \leq k-1-\des(\pi).\] The bijection between these sets of solutions is given by $i'_s = i_s - a(s)$ where $a(s)$ is the number of descents to the left of $s$. Therefore the $q$-order polynomial is given by
\begin{eqnarray*}
\Omega_{\pi}(q;k) & = & \sum_{\substack{0\leq i_1 \leq \cdots \leq i_n \leq k-1 \\ s\in \Des(\pi) \Rightarrow i_s < i_{s+1}}} \left(\prod_{s=1}^{n} q^{i_s}\right)\\
 & = & \sum_{0\leq i'_1 \leq \cdots \leq i'_n \leq k-1-\des(\pi)} \left(\prod_{s=1}^{n} q^{i'_s + a(s)}\right)\\
 & = & q^{\sum_{s=1}^{n} a(s)}\cdot\left(\sum_{0\leq i'_1 \leq \cdots \leq i'_n \leq k-1-\des(\pi)} \left(\prod_{s=1}^{n} q^{i'_s}\right)\right).
\end{eqnarray*}
The sum of all $a(s)$ can be expressed as $\displaystyle\sum_{s\in \Des(\pi)}(n-s)$, which is sometimes referred to as the \emph{comajor index}, denoted $\comaj(\pi)$.\footnote{The \emph{major} index of a permutation is $\displaystyle\sum_{s\in \Des(\pi)}\!\!\!\!s$. Indeed, had we made our $P$-partitions order reversing, as in Stanley's original definition of a $P$-partition, we would have gotten $q^{\maj(\pi)}$ rather than $\comaj$ above.} The rest of the sum is now recognizable as a $q$-binomial coefficient. In summary, we have
\begin{equation}\label{eq:qorder}
\Omega_{\pi}(q;k) = q^{\comaj(\pi)}\binom{k+n-1-\des(\pi)}{n}_q .
\end{equation}

Now we will give a $q$-analog of Theorem \ref{thm:ges}. We would like to have a $q$-analog for Theorem \ref{thmcyc} as well, but if such a theorem exists, it does not seem to be immediate given the method of proof in this paper. Recall that the proof of Theorem \ref{thmcyc} hinged on the observation that cyclically permuting a permutation $\pi \in \mathfrak{S}_n$ leaves the number of cyclic descents unchanged: $\cdes(\pi) = \cdes(\pi \omega^i)$ for any $i = 0, 1, \ldots, n-1$. Unfortunately, cyclically permuting $\pi$ has a more subtle effect on the comajor index.

Moving on, we define the $q$-version of the structure polynomial, \[\phi(q;x) = \sum_{\pi\in\mathfrak{S}_n} q^{\comaj(\pi)}\binom{x+n-1-\des(\pi)}{n}_q \pi .\]

\begin{figure} [t]
\centering
\includegraphics[scale=.7]{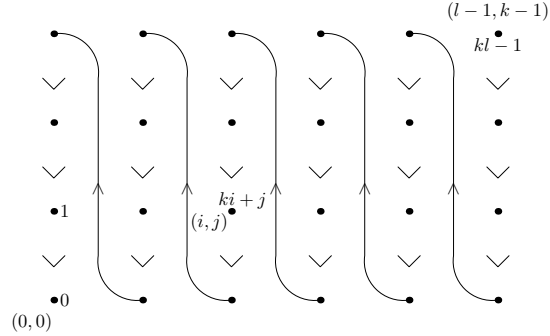}
\caption{Weights in the lexicographic order.}
\end{figure}

\begin{thm}\label{thm:qsymmetric}
As polynomials in $x$ and $y$ (and $q$) with coefficients in the group algebra we have \[ \phi(q;x) \phi(q^{x};y) = \phi(q;xy).\]
\end{thm}

\begin{proof}
The proof will follow nearly identical lines of reasoning as in the ordinary ($q =1$) case. Here we sketch the proof with emphasis on the major differences. Again, we will decompose the coefficient of $\pi$: \[q^{\comaj(\pi)}\binom{kl+n-1-\des(\pi)}{n}_q\] By \eqref{eq:qorder}, we have that the coefficient of $\pi$ is the order polynomial $\Omega_{\pi}(q;kl)$ so we will examine the $\pi$-partitions $f: [n]\to \{0,1,\ldots,l-1\}\times \{0,1,\ldots,k-1\}$. Notice that we are still mapping into a set with $kl$ elements. As before we impose the lexicographic ordering on this image set. To ensure that we keep the proper powers of $q$, we think of the order polynomial now as:\[\Omega_{\pi}(q;kl) =\sum_{\substack{(0,0)\leq (i_1,j_1) \leq \cdots \leq (i_n,j_n) \leq (l-1,k-1) \\ s\in \Des(\pi) \Rightarrow (i_s,j_s) < (i_{s+1},j_{s+1})}} \left(\prod_{s=1}^{n} q^{ki_s+j_s}\right).\] We have given each point $(i,j)$ the weight $ki+j$ so that the weight corresponds to the position of the point in the lexicographic ordering on $\{0,1,\ldots,l-1\}\times\{0,1,\ldots,k-1\}$. We now proceed exactly as in the proof of Theorem \ref{thm:ges}.
\begin{eqnarray*}
\Omega_{\pi}(q;kl) &=& \sum_{\substack{(0,0)\leq (i_1,j_1) \leq \cdots \leq (i_n,j_n) \leq (l-1,k-1) \\ s\in \Des(\pi) \Rightarrow (i_s,j_s) < (i_{s+1},j_{s+1})}} \left(\prod_{s=1}^{n} q^{ki_s+j_s}\right)\\
 & = & \sum_{I\subset[n-1]}\left(\sum_{\substack{ 0\leq i_1 \leq \cdots \leq i_n \leq l-1\\ s\in I \Rightarrow i_s < i_s+1}} q^{ki_s} \right) \left(\sum_{\sigma \in \mathcal{L}(P_I)} \Omega_{\sigma}(q;k) \right)\\
& = & \sum_{\sigma\tau = \pi} \Omega_{\sigma}(q;k) \Omega_{\tau}(q^{k};l),
\end{eqnarray*}
as desired.

\end{proof}

In order to give a $q$-analog of Theorem \ref{thm:cho}, we define the $q$-order polynomial for a signed permutation $\pi \in \mathfrak{B}_n$ as
\begin{eqnarray*}\Omega_{\pi}(q;k) & = & \sum_{\substack{0\leq i_1 \leq \cdots \leq i_n \leq k \\ s\in \Des(\pi) \Rightarrow i_s < i_{s+1}}} \left( \prod_{s=1}^{n} q^{i_s} \right), \\
\mbox{which as before can be rewritten,} & & \\
 & =  & \sum_{0\leq i'_1 \leq \cdots \leq i'_n \leq k-\des(\pi)} \left(\prod_{s=1}^{n} q^{i'_s + a(s)}\right)\\
 & = & q^{\sum_{s=1}^{n} a(s)}\cdot\left(\sum_{0\leq i'_1 \leq \cdots \leq i'_n \leq k-\des(\pi)} \left(\prod_{s=1}^{n} q^{i'_s}\right)\right) \\
 & = & q^{\comaj(\pi)}\binom{k+n-\des(\pi)}{n}_{q}.
\end{eqnarray*}
Now we define the type B $q$-structure polynomial, \[\phi(q;x) = \sum_{\pi\in\mathfrak{B}_n} q^{\comaj(\pi)} \binom{x+n-\des(\pi)}{n}_q \pi.\]
\begin{thm}\label{thm:qhyper}
The following relation holds as polynomials in $x$ and $y$ (and $q$) with coefficients in the group algebra of the hyperoctahedral group: \[ \phi(q; x)\phi(q^{2x+1}; y)= \phi(q; 2xy+x+y),\] or upon substituting \[x \leftarrow (x-1)/2,\] \[y \leftarrow (y-1)/2,\]
then \[ \phi(q;(x-1)/2)\phi(q^{x};(y-1)/2) =\phi(q;(xy-1)/2).\]
\end{thm}
\begin{proof}
We will omit most of the details, but the crucial step is to keep the proper exponent on $q$. We give to each point $(i,j)$ the weight $(2k+1)i + j$ so that the weight corresponds to the position of the point in the lexicographic order on the set $\{0,1,\ldots,l\}\times \{-k,\ldots ,-1,0,1,\ldots,k\}$. The proof is outlined in two steps below. For any $\pi$ and any pair of positive integers $k,l$,
\begin{eqnarray*}
\Omega_{\pi}(q; 2kl+k+l) & = & \sum_{\substack{ (0,0)\leq (i_1,j_1) \leq \cdots \leq (i_n,j_n) \leq (l,k) \\ s \in \Des(\pi) \Rightarrow (i_s,j_s) < (i_{s+1},j_{s+1}) }} \left( \prod_{s=1}^{n} q^{(2k+1)i_s + j_s} \right) \\
 & = & \sum_{\sigma\tau = \pi} \Omega_{\sigma}(q; k) \Omega_{\tau}(q^{2k+1};l).
\end{eqnarray*}
\end{proof}

The augmented version of the $q$-order polynomial gives us $q$-analogs for Theorem \ref{conj1} and Theorem \ref{thm:hypermod}. For a signed permutation $\pi \in \mathfrak{B}_n$ we have
\[
\Omega^{(a)}_{\pi}(q;k) = \sum_{\substack{0 \leq i_1 \leq \cdots \leq i_n \leq k\\ s \in \aDes(\pi) \Rightarrow i_s < i_{s+1} } } \left( \prod_{s = 1}^{n} q^{i_s} \right) = q^{\acomaj(\pi)} \binom{k + n - \ades(\pi)}{n}_q,
\]
where if $a(s)$ is the number of descents of $\pi$ to the left of $s$, then the \emph{augmented comajor index}, $\acomaj(\pi)$, is the sum over all $s$ of the numbers $a(s)$. Define \[\psi(q;x) = \sum_{\pi \in \mathfrak{B}_n} q^{\acomaj(\pi)} \binom{x + n - \ades(\pi)}{n}_q \pi.\]

\begin{thm}\label{thm:qaug}
As polynomials in $x$ and $y$ (and $q$) with coefficients in the group
algebra of the hyperoctahedral group we have
\[\psi(q;x)\psi(q^{2x};y) = \psi(q;2xy),\] or upon substituting \[x \leftarrow x/2,\] \[y \leftarrow y/2,\]
then \[\psi(q;x/2)\psi(q^x; y/2) = \psi(q; xy/2).\]
\end{thm}

\begin{figure} [h]
\centering
\includegraphics[scale=.7]{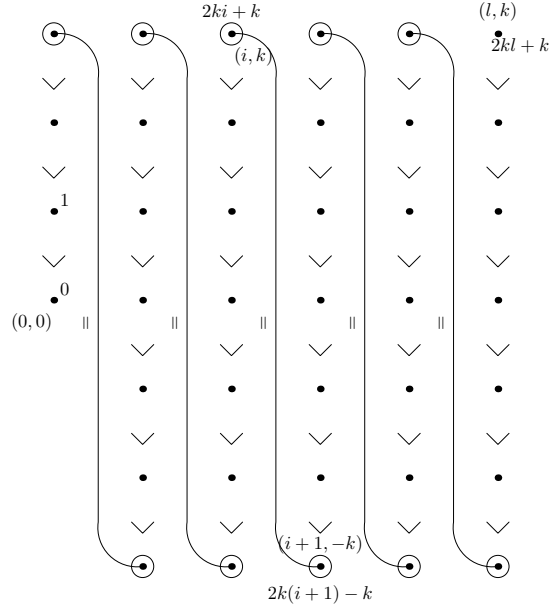}
\caption{Weights in the augmented lexicographic order.}
\end{figure}

\begin{proof}

Here the key is to give the integer pairs $(i,j)$ the proper weight in the augmented lexicographic ordering. If we take $2ki + j$ as the weight of the point $(i,j)$ then we get that the points $(i,k)$ and $(i+1,-k)$ have the same weight. So the weight does indeed correspond to the position of $(i,j)$ in the augmented lexicographic ordering. Everything else follows as in the proof of Theorem \ref{conj1}. For any $\pi \in \mathfrak{B}_n$ and any pair of positive integers $k,l$,
\begin{eqnarray*}
\Omega^{(a)}_{\pi}(q; 2kl) & = & \sum_{\substack{ (0,0)\leq (i_1,j_1) \leq \cdots \leq (i_n,j_n) \leq (l,0) \\ s \in \aDes(\pi) \Rightarrow (i_s,j_s) < (i_{s+1},j_{s+1}) }} \left( \prod_{s=1}^{n} q^{2k i_s + j_s} \right) \\
 & = & \sum_{\sigma\tau = \pi} \Omega^{(a)}_{\sigma}(q; k) \Omega^{(a)}_{\tau}(q^{2k};l).
\end{eqnarray*}

\end{proof}

\begin{thm}\label{thm:qhypermod}
As polynomials in $x$ and $y$ (and $q$) with coefficients in the group algebra of the hyperoctahedral group we have
\[\psi(q;x)\phi(q^{2x};y) = \psi(q;2xy+x),\] upon substituting \[x \leftarrow x/2,\] \[y \leftarrow (y-1)/2,\]
then \[\psi(q;x/2)\phi(q^x;(y-1)/2) = \psi(q;xy/2).\]
\end{thm}

\begin{proof}
Because we exploit the augmented lexicographic order in the proof of Theorem \ref{thm:hypermod} (the $q=1$ case), we will use the same weighting scheme as in the proof of Theorem \ref{thm:qaug} for the points $(i,j)$. We have for any $\pi \in \mathfrak{B}_n$ and any pair of positive integers $k,l$,
\begin{eqnarray*}
\Omega^{(a)}_{\pi}(q; 2kl+k) & = & \sum_{\substack{ (0,0)\leq (i_1,j_1) \leq \cdots \leq (i_n,j_n) \leq (l,k) \\ s \in \Des(\pi) \Rightarrow (i_s,j_s) < (i_{s+1},j_{s+1}) }} \left( \prod_{s=1}^{n} q^{2k i_s + j_s} \right) \\
 & = & \sum_{\sigma\tau = \pi} \Omega^{(a)}_{\sigma}(q; k) \Omega_{\tau}(q^{2k};l).
\end{eqnarray*}
\end{proof}

\section{Future work}\label{sec:fut}

We hope to be able to apply the method of $P$-partitions to study descent algebras more generally. In his thesis work, Vic Reiner \cite{R2} defined $P$-partitions purely in terms of the root system of a Coxeter group. His definition may provide a way to obtain a theorem that generalizes Theorems \ref{thm:ges} and \ref{thm:cho} to any Coxeter group.

Another avenue that we have already begun to investigate is the study of \emph{peak algebras}. A \emph{peak} of a permutation $\pi \in \mathfrak{S}_n$ is a position $i$ such that $\pi(i-1) < \pi(i) > \pi(i+1)$. Kathryn Nyman \cite{N} showed that the span of sums of permutations with the same set of peaks forms a subalgebra of the group algebra. Later, Aguiar, Bergeron, and Nyman \cite{AguBerNym} explored many properties of this and related peak algebras for the symmetric group, including connections to descent algebras of types B and D. A key tool in \cite{N} was John Stembridge's \emph{enriched} $P$-partitions \cite{Ste}. In forthcoming work, we use Stembridge's tool to study subalgebras formed by the span of sums of permutations with the same \emph{number} of peaks, obtaining results similar to the main theorems of this paper.

A meta-question for this area of research is: what is so special about descents (or peaks)? Why should this way of grouping permutations give rise to a subalgebra at all? In terms of permutation patterns, a descent is an occurrence of the pattern 21, and a peak is an occurrence of the pattern 132 or 231. What distinguishes descents and peaks from, say, grouping permutations according to instances of 1324 (which doesn't give a subalgebra)?

\end{document}